\documentclass[12pt]{article}

\usepackage{amsthm}
\usepackage{amsmath,amssymb}
\usepackage{verbatim}
\usepackage{xcolor}
\usepackage{graphicx,wrapfig}
\usepackage[affil-it]{authblk}

\usepackage{caption}

\newtheorem{theorem}{Theorem}[section]

\newtheorem{lemma}[theorem]{Lemma}

\newtheorem{remark}[theorem]{Remark}
\newtheorem{example}[theorem]{Example}
\numberwithin{equation}{section}

\def\be{\begin{equation}}
\def\ee{\end{equation}}
\def\bes{\begin{equation*}}
\def\ees{\end{equation*}}

\definecolor{dgreen}{rgb}{0, 0.6, 0.1}

\newcommand{\supp}{\mathop{{\rm supp}}}

\renewcommand{\P}{\mathbb{P}}
\newcommand{\E}{\mathbb{E}}
\newcommand{\Rd}{\mathbb{R}^d}
\newcommand{\R}{\mathbb{R}}

\newcommand{\barr}{\begin{array}{rcl}}
\newcommand{\earr}{\end{array}}

\newcommand{\bR}{\mathbb{R}}
\newcommand{\bZ}{\mathbb{Z}}
\newcommand{\bH}{\mathbb{H}}
\newcommand{\bP}{\mathbb{P}}

\newcommand{\pd}{\partial}
\newcommand{\ol}{\overline}
\newcommand{\al}{\alpha}

\begin{document}

\title{Some Boundary Harnack Principles With Uniform Constants}


\author{Martin T. BARLOW%
	\thanks{Research partially supported by NSERC (Canada). E-mail: \texttt{barlow@math.ubc.ca}}}
\affil{Department of Mathematics, \\University of British Columbia, \\Vancouver B.C., Canada V6T 1Z2.}

\author{Deniz KARLI%
	\thanks{Research partially supported by NSERC (Canada) and BAP 20A101 Grant of I\c{s}\i k University (Turkey). Email: \texttt{deniz.karli@gmail.com}}}
\affil{Department of Mathematics, \\I\c{s}\i k University, \\\c{S}ile, Istanbul, Turkey 34980.}


\date{}

\maketitle


\begin{abstract}
	We prove two versions of a boundary Harnack principle in which 
	the constants do not depend on the domain by using probabilistic methods.\\
	
	{\it Keywords: Boundary Harnack principle; Harmonic functions; Brownian motion; Harnack inequality }
	
	{MSC 60J45 \and MSC  31C05 \and MSC 42A61}
\end{abstract}

\section{Intoduction}\label{intro}

A boundary Harnack principle (BHP) gives a result of the following general type.
Let $D$ be a domain in $\bR^d$, and $\xi \in \pd D$, satisfying suitable 
properties.
Let $r>0$, $a_0 \ge 2$,

$B_1 = B(\xi, a_0 r)$ and $B_2 = B(\xi,r)$; here $B(.,.)$ denote
the usual Euclidean balls. Then there exists a constant $C_{D}$ such that
if $u,v$ are positive harmonic functions on $B_1 \cap D$
vanishing on $\pd D \cap B_1$, one has
\be \label{e:basicBHP}
\frac{ u(x)/v(x) }{ u(y)/v(y)} \le C_D \hbox { for } x,y \in D \cap B_2. 
\ee
A BHP of this kind is called in \cite{Aik01} a {\em uniform} BHP, and in 
\cite{L-SC} a {\em scale invariant} BHP. Here `uniform' or 
`scale invariant' refers to the fact that the constant $C_D$ does not depend on $r$.
For Lipschitz  domains $D$ the scale invariant BHP was proved independently by
Ancona, Dahlberg and Wu in 
\cite{Anc78,Dahl,Wu}. This was extended to NTA domains by 
Jerison and Kenig \cite{JK}. Bass, Burdzy and Banuelos \cite{BB,BBB} used
probabilistic methods to obtain a 
BHP for H\"older domains, but their BHP is not uniform. In \cite{Aik01}
a scale invariant BHP is proved for uniform domains in $\bR^d$, and in
\cite{aikawa1} this is extended to John domains.
See the papers \cite{Aik01,AikC,L-SC} for a further discussion on the history of 
the BHP, and the various different kinds of BHP.

In the above `harmonic  function' refers to functions which are harmonic with respect
to the usual Laplacian operator in $\bR^d$. (These functions are
harmonic with respect to the infinitesimal generator of the semigroup of
standard Brownian motion in $\bR^d$.) 
Recent papers have studied functions which are harmonic with respect to the 
generators of more general diffusion processes -- see \cite{L-SC,L14} and the 
references therein.

In all these results the constant $C_D$ depends on the domain $D$.
For the standard Laplacian it
is clear that such dependence is necessary, since the BHP does not hold for
all domains $D \subset \bR^d$. 
(See however \cite{bkk} where a  BHP  with constants independent
of the domain is proved for harmonic functions with respect to fractional Laplacians.)

This paper originates in the work of Masson \cite{Mas}, where a boundary estimate
with a constant $C_D$ {\em not depending on $D$} was needed --
see \cite[Proposition 3.5]{Mas}.
Masson's work was in the context of discrete potential theory for $\bZ^2$.
Let $S^x=(S^x_k, k \ge 0)$ be the simple random walk on $\bZ^2$, started
at $x$, and write $S=S^0$. 
Write $\bZ^2_- =\{ (x_1, x_2) \in \bZ^2: x_1 \le 0 \}$, and let 
$Q(x,n)=\{ y \in \bZ^2: |x-y | \le n\}$. Let $N \ge 1$, 
$ K \subset Q(0,N) \cap \bZ^2_-$, and $D = Q(0,N)-K$. (The case of interest is 
when $0 \in K$.) 
Let $\tau^+=\min\{ k \ge 1: S_k \not\in D\}$, and
$F =\{ S_{\tau^+} \in \bZ^2 - K \}$, so that $F$ is the event that 
$S$ leaves $Q(0,N)$ before hitting $K$.
Let $W =\{ (x_1, x_2):  0 \le |x_2| \le x_1 \}$, so that $W$ is a cone
with vertex $(0,0)$ and angle $\pi/4$. Massons's theorem is that there
exists $p_0>0$, independent of $N$ and $K$, such that $\bP( S_{\tau^+} \in W | F) \ge p_0$.
This result extends to give also
\be \label{e:maslb}
\bP( S^x_{\tau^+} \in W | F) \ge p_1 \hbox { for $x=(x_1,x_2) \in Q(0,N/16)$ with $x_1\ge 0$ }.
\ee
The fact that the constants $p_i$ do not depend on the structure of $K$ is essential in the
context of \cite{Mas}, since $K$ is a random path (actually a loop erased random walk), 
and the estimate \eqref{e:maslb} was needed for all possible $K$.

Although the connection with BHP is not made in \cite{Mas}, this
result is clearly of BHP type. For $x \in \bZ^2$ let 
$\tau=\min\{ k \ge 0: S_k \not\in D\}$, and define the functions
$$  v(x) = \bP( S^x_\tau \in K^c ), \quad u(x) = \bP( S^x_\tau \in K^c \cap W  ). $$ 
These are (discrete) harmonic in $D$, and
$\bP( S^x_\tau \in W | F) = {u(x)}/{v(x)}$. 
Since $u \le v$ it is immediate from \eqref{e:maslb} that 
\be \label{e:masbhp}
\frac{ u(x)/v(x) }{u(y)/v(y)} \le p_1^{-1},
\text{ for } x,y \in Q(0, N/16) \cap (\bZ^2 - \bZ^2_-).
\ee

Thus we have a BHP for the specific functions $u,v$
in which the constant $C_D=p_1^{-1}$
does not depend on $K$; the price is that the inequality only holds
for those $x \in Q(0,N/16)$ with $x_1> 0$.

{\centering \includegraphics[width=0.8\columnwidth]{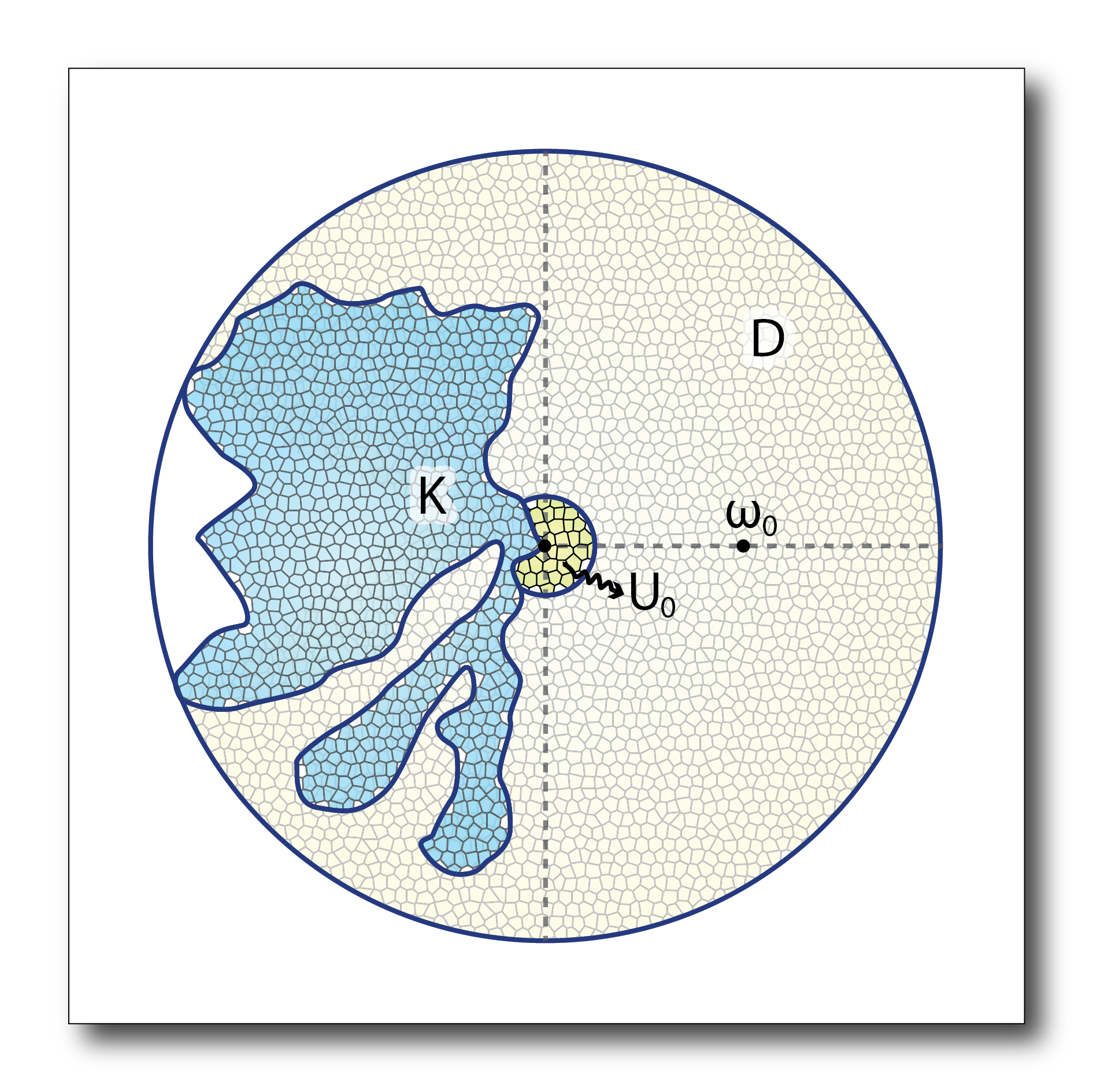}
	\captionof{figure}{The sets $K$ and $U_0$. }
}

\medskip
Our first result is a BHP in two dimensions which holds with a constant independent 
of the domain. 
We write $B(x,r)=\{ y: |x-y|<r\}$ for Euclidean balls with center at $x$ and radius $r$,
and set 
$$ \bH_- =\{ (x_1, \dots, x_d) \in \bR^d : x_1 < 0 \}, $$
for the open left-halfspace, and define the  open right-halfspace $\bH_+$ analogously. 

\begin{theorem} \label{T:ubhp2} 
	Let $d=2$, let $K\subset B(0,1) \cap \ol \bH_-$ be connected and relatively closed in 
	$B(0,1)$, and let $D=B(0,1) - K$. Let $U_0$ be the connected component
	of $B(0, \tfrac{1}{16})\cap D$ which contains $(\tfrac{1}{32},0)$. 
	There is a positive constant $C_0$, independent of $K$, such that
	if $u$ and $v$ are positive and harmonic on $D$, are continuous on $\overline{D}$, 
	and vanish on $K$ then
	\be \label{e:bhp2b}
	\frac{u(x)/v(x)}{u(y)/v(y)} < C_0\quad \text{ for } x,y \in U_0. 
	\ee
\end{theorem}

\smallskip A key estimate for the proof is the Carleson estimate Lemma \ref{rulce},
which is proved by a path-crossing argument. This estimate relies on the fact that
$K$ is connected, and does not generalize to $d\ge 3$. Examples at the end of Section 
\ref{S:BHPnl} show that Theorem \ref{T:ubhp2} does not hold in general if $d \ge 3$, or 
if $K$ is not connected, and that the inequality \eqref{e:bhp2b} cannot be extended to
$x,y \in B(0, \tfrac{1}{16})\cap D$.

In Section \ref{S:BHPuv} we extend Masson's result to Euclidean space, and
prove it for $d \ge 2$. (The result is trivial for $d=1$).

\smallskip
Throughout this paper we write 
$X=(X_t, t \in [0,\infty), \bP^x, x \in \bR^d)$ for Brownian motion in $\bR^d$;
$\bP^x$ is the law of $X$ started at $x$. For a set $A\subset \bR^d$
we define
$$  T_A = \inf\{ t\ge 0: X_t \in A \}, \quad \tau_A = T_{A^c} 
= \inf\{ t\ge 0: X_t \not\in A \}. $$
We will use the classic and probabilistic definitions of harmonic functions interchangeably. 
We call a function $h$ {\em harmonic} in a domain $A$ if $h$ is locally integrable and for all 
$x\in A$ and all $r<{\rm dist}(x,\partial A)$, 
$$h(x)=\frac{1}{|B(0,r)|}\int_{B(x,r)}h(y)dy.$$
Equivalently, $h$ is harmonic in $A$ if $h(X_{t\wedge \tau_A})$ is a martingale.

When we use notation such as $C=C(\alpha,d)$ this will mean that the
(positive) constant $C$ depends only on the parameters $\alpha$ and $d$.\\

\section{BHP for positive harmonic functions in $d=2$}
\label{S:BHPnl}

In this section we prove Theorem \ref{T:ubhp2}. We begin with a general lemma
concerning harmonic functions in a bounded domain in $\bR^d$.

\begin{lemma}\label{path_connected}
	Let $D$ be a bounded path connected domain in $\R^d$, with $d\ge 1$, and 
	$f$ be a non-negative harmonic function in $D$ which is continuous on $\overline{D}$. 
	Let $c_0>0$. If there exists $x_0\in D$ such that $f(x_0)>c_0$ then the set 
	$$\partial D_0:=\{x\in \partial D: f(x)>c_0\}\not=\emptyset.$$
	Moreover, there is a  path $\gamma$ from $x_0$ to a 
	point  $x^* \in \partial D_0$ such that 
	$\gamma\setminus \{x^*\}  \subseteq D$ and $f(x)>c_0$ for all $x\in \gamma$.
\end{lemma}

\begin{proof}
	Since $f$ is a non-negative harmonic function in $D$ and it is continuous on $\ol D$, 
	by the Maximum Principle for harmonic functions, there exists a point $x_1\in \partial D$ such that 
	$f(x_1)\geq f(x_0)>c_0.$ Hence $\partial D_0\not=\emptyset$.
	Define $\Gamma$ to be the collection of all paths $\gamma$ from $x_0$ to a point in $\partial D_0$ 
	so that $\gamma\subseteq D$ and $f(y) > c_0$ for all $y \in \gamma$.
	We will show that $\Gamma\not=\emptyset$. 
	Denote Brownian motion starting at the point $x_0$ by $(X_t,\P^{x_0})$, and stopping times
	$$ T_0=\inf\{ t\ge 0 :f(X_t)\leq c_0\}\quad \mbox{and}\quad T=T_0\wedge \tau_D$$
	where $\tau_D$ is the first exit time of $X_t$ from the domain $D$. Then
	$$ c_0 < f(x_0)=\E^{x_0}(f(X_T))=c_0\,\P^{x_0}(T_0\leq \tau_D)+\E^{x_0}(f(X_{\tau_D});\tau_D<T_0).$$
	If $\P^{x_0}(T_0\leq \tau_D)=1$ then $f(x_0)=c_0$ which contradicts the assumption. 
	Hence $\P^{x_0}(T_0\leq \tau_D)<1$ and so $\P^{x_0}(T_0> \tau_D)>0$. 
	Notice that on the event $\{T_0> \tau_D\}$ the path of the Brownian motion is in $\Gamma$;
	since this event has non-zero probability, $\Gamma\not=\emptyset$. 
\end{proof}

Now let $K \subset B(0,1)$ satisfy the hypotheses of the Theorem.
First, we will create a domain enclosing the region $U_0$ inside the domain $D$. 
For this purpose, we let $r_1=1/4$ and $r_2=3/4$. 
We begin by assuming that $$K \cap \pd B(0,r_j) \neq \emptyset$$ for $j=1,2$. (See Fig. \ref{fig2}.)
We also write $r_0={1}/{16}$; this is the radius of the smaller balls used in Carleson estimate below. 

{  \centering
	\includegraphics[width=0.8\columnwidth]{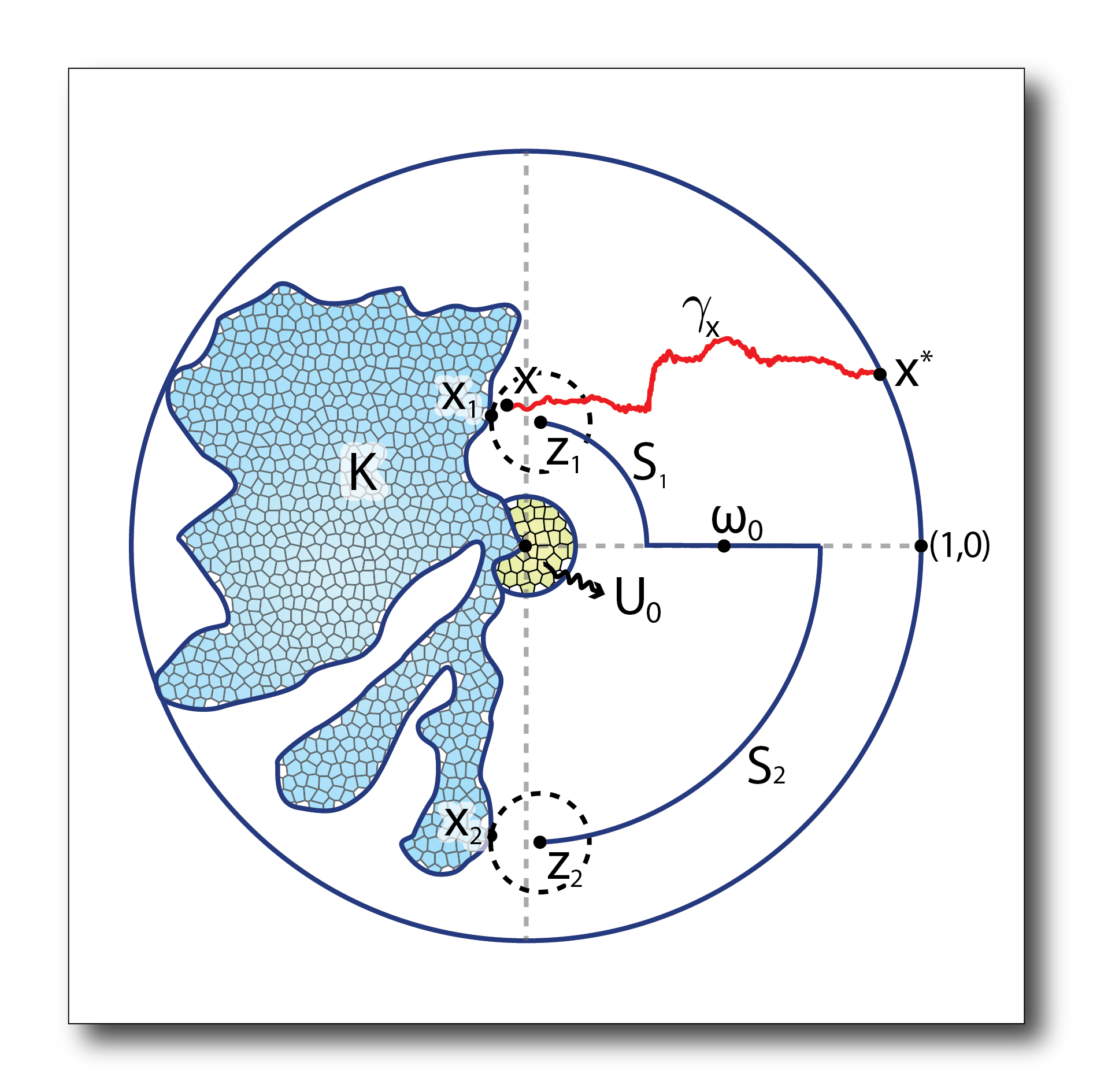}
	\captionof{figure}{The points $x_j$, $z_j$ and curves $S_j$ }\label{fig2}
}

\medskip
We write $w_0 = (\tfrac12,0)$ and define two hitting angles $\theta_1$ and $\theta_1$ as follows:
\begin{align*}
	\theta_1& =\inf\{\theta\in(0,2\pi): \overline{B(r_1e^{i\theta},r_0)}\cap K \not=\emptyset\},\\
	\theta_2 & =\sup\{\theta\in(-2\pi,0): \overline{B(r_2 e^{i\theta},r_0)}\cap K \not=\emptyset\}.
\end{align*}

\medskip
Write $z_j = r_je^{i \theta_j}$ and let $x_j\in K\cap\overline{B(z_j, r_0 )}$ for $j=1,2$. 
Note that the balls $B(x_j, r_0)$, $j=1,2$ are disjoint and the distance between the points 
$x_j$ and $z_j$ is $r_0$, 
Let $x_j'$ be the midpoint on the line segment $[x_j, z_j]$,
$L_j$  be the the line segment $[x'_j, z_j]$, and $x''_j$ be the midpoint of $L_j$.
We define the annulus $A_j = B(x_j, r_0) -  \ol{B(x_j, \frac{r_0}{2})} $.		
Let $U_j$ 	be the connected component of 
$ B(x_j, \frac{r_0}{2}) - K$ which contains the open line segment between 
$x_j$ and $x'_j$. As $U_j$ is an open connected subset of $\bR^2$ it is also path-connected.



We begin by proving two Lemmas; the first ensures existence of certain paths in the domain,
while the second is a rooted local Carleson estimate. 
Recall from the statement of Theorem \ref{T:ubhp2} that $U_0$ is the connected 
component of $B(0, \tfrac{1}{16})\cap D$ which contains $(\tfrac{1}{32},0)$.

\begin{lemma} [Carleson estimate]  \label{rulce}
	Let $K$, $D$, $U_1$, $U_2$ be as above. 
	Let $z \in U_0$ and $u$ be as in Theorem {\ref{T:ubhp2}}.
	There exists a constant $C$, independent of $K$, $z$ and $u$, such that
	$$u(y)\leq C u(w_0) \quad \mbox{and} \quad G_D(z, y)\leq C\, G_D(z,w_0),  \quad y \in U_j .$$ 
\end{lemma}

\def\gam{\gamma}

\begin{proof} We will prove this for $U_1$; the same argument also applies to $U_2$. 
	
	The set $\pd B(z_1 , r_0) \cap A_1 $
	consists of the union of two connected arcs; denote these 
	$\gamma_2$ and $\gamma_3$, labelled so that going anticlockwise round $A_1$ we meet 
	the arcs $L_1, \gamma_2, \gamma_3$ in order. (See Fig. \ref{fig3}.)
	
	{\centering
		\includegraphics[width=0.8\columnwidth]{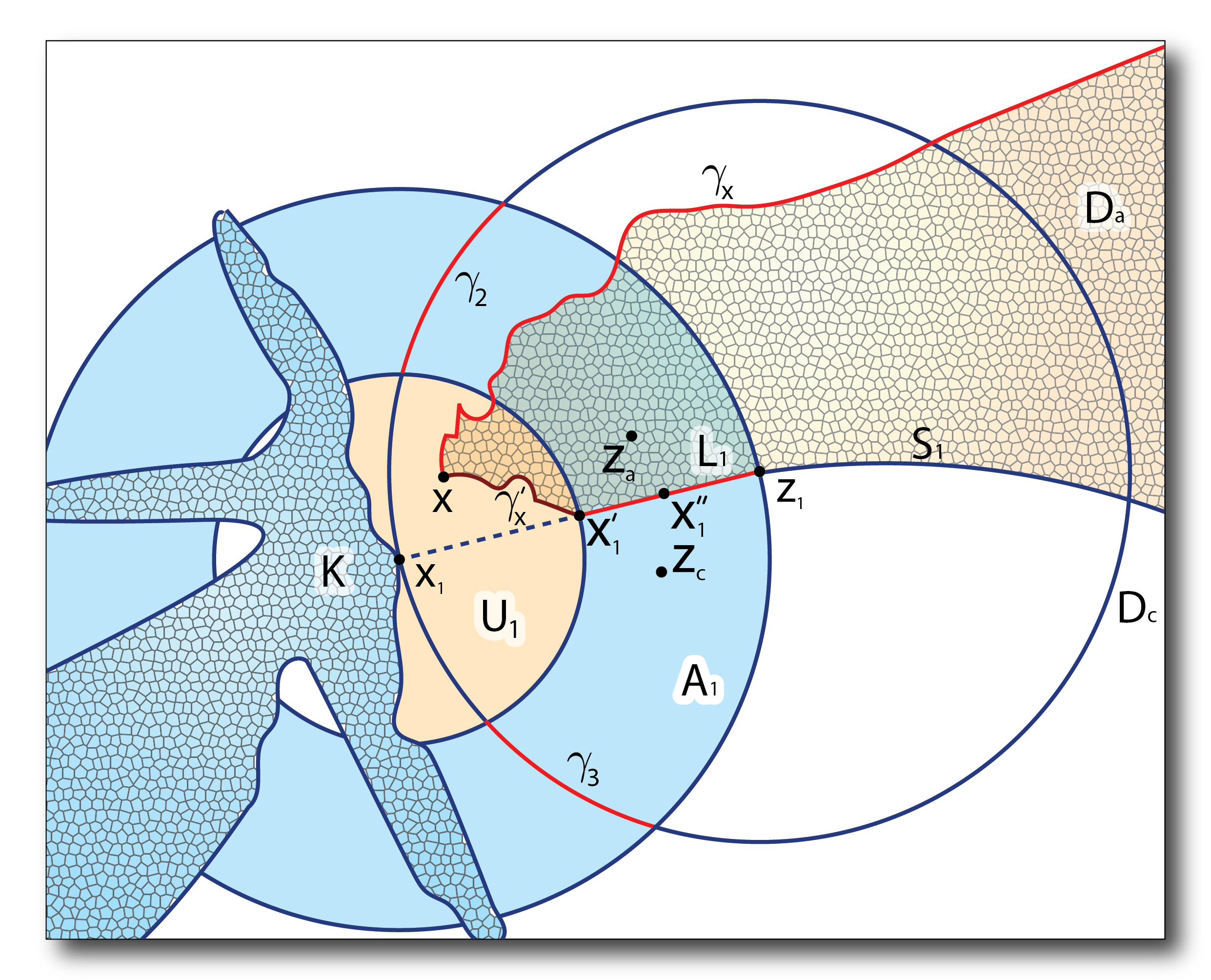}
		\captionof{figure}{ We zoom into the region around the points $x_1$ and $z_1$. 
			A similar region exists around $x_2$ and $z_2$.}\label{fig3}
	}
	
	\medskip
	Let $H_2$ be the event that a Brownian motion $X$, started on the line $L_1$
	stays inside $A_1$ until after it has hit in order the sets $\gamma_2$, $\gamma_3$, and then $L_1$. 
	More precisely if we set
	\begin{align*}
		T_{23} &= \inf\{ t \ge 0: X_t \in \gamma_2 \cup \gamma_3 \}, \\
		T_{31} &= \inf\{ t \ge T_{23} : X_t \in \gamma_3 \cup L_1 \}, \\
		T_{12} &= \inf\{ t \ge T_{31} : X_t \in L_1 \cup \gamma_2 \}, 
	\end{align*}
	then
	$$ H_2 =\{ T_{23} < T_{31} < T_{12} < \tau_{A_1}, 
	X_{T_{23}} \in \gamma_2, X_{T_{31}} \in \gamma_3, X_{T_{12}} \in L_1 \}. $$
	Let $H_3$ be the similar event with the roles of $\gamma_2$ and $\gamma_3$ interchanged.
	As the sets $L_1, \gamma_2, \gamma_3$ are separated by a distance $c r_0$, and using the symmetry 
	of the set, there exists $p_1>0$ such that
	$$ \bP^{x_1''}( H_2)  = \bP^{x_1''}( H_3) = p_1. $$
	By the Harnack inequality there exists a constant $C_2$ such that if $h$ is non-negative and harmonic in 
	$B(z_1, r_0)$ then 
	\be \label{e:harm2}
	h(y) \le  C_2 \, h(x''_1) \hbox{ for all } y \in B(x''_1, \frac{r_0}{3}). 
	\ee

	Now let $C_1 = \max\{ 2/p_1, C_2\}$.
	We consider first the case when $f=u$. It is enough to prove
	\begin{equation}
		f(x) \le C_1 f( x_1''), \quad \hbox{ for } x \in U_1,
	\end{equation} 
	since then by using the Harnack inequality in a chain of balls on the arc 
	$\{ r_1 e^{i\theta}, 0 \le \theta \le \theta_1\}$ and the line $\{ (t,0), r_1\le t \le \tfrac12 \}$
	we have $f(x''_1) \le c f(w_0)$.
	
	If $f(x) \le C_1 f(x''_1)$ for all $x \in U_1$ then we are done.
	So suppose there exists $x \in  U_1$ with  $f(x) > C_1 f(x''_1)$. 
	As $f$ is harmonic and non-negative in $D$, 
	by Lemma \ref{path_connected}, there exists a path $\gamma_x$ from $x$ to  a point $x^* \in \pd B(0,1)$  such that
	$f(y) >   C_1 f(x''_1)$ for all $y \in \gamma_x$.

	Suppose that there exists $y \in \gamma_x \cap B(x''_1, \frac{r_0}{3})$.
	Then using \eqref{e:harm2} we have
	$f(x''_1) \ge C_2^{-1} f(y) > C_2^{-1} C_1  f(x''_1)>  f(x''_1)$, a contradiction. 
	So we have that $\gamma_x \cap B(x''_1, \frac{r_0}{3}) = \emptyset$.

	We now define a path $\gamma$ in $\ol D$ between $x^*$ and $(1,0)$ as follows,
	which includes the points $x^*, x, x'_1, z_1, (r_1,0), (1,0)$.
	By the definition of $U_1$, there is a path $\gamma'_x$ in $B(x_1, \frac{r_0}{2}) - K$
	connecting $x$ and $x_1'$. Let $S_1$ be the path $(r_1 e^{i \theta}, 0 \le \theta \le \theta_1)$
	which connects $(r_1,0)$ and $z_1$, and $S_2$ be the line segment between $(r_1,0)$ and $(1,0)$. (See Fig. \ref{fig2} and Fig. \ref{fig3}.)
	Then $\gamma$ consists of the  concatenation (with appropriate orientations) of $\gamma_x, \gamma'_x, L_1, S_1, S_2$.
	(The path $\gamma$ is not necessarily a simple curve -- it may have multiple points.)
	We also write $\gamma$ for the set of points in this path. By the construction of $\gamma$ we have 
	that $\gamma \cap K = \emptyset$. 
	
	Let $z_a$ and $z_c$ be two points close to
	$x''_1$ on opposite side of the line segment $L_1$, in the anticlockwise and clockwise directions
	respectively, and let $D_a$ and $D_c$ be the connected 
	components of $B(0,1) - \gamma$ which contain $z_a$ and $z_c$ respectively. 
	Since the path $\gamma$ inside $B(x_1'', r_0/3)\cap A_1$ just consists of the line segment $L_1$ without 
	its endpoints,
	the components $D_a$ and $D_c$ are distinct. 
	(See Remark \ref{R:winding} below for more details.)
	Hence, as $K \cap \gamma = \emptyset$ and $K$ is connected, 
	at most one of $K \cap D_a$, $K \cap D_c$ is non-empty. 
	
	We suppose that
	$K \cap D_a = \emptyset$. 
	By the construction of $\gamma$, if the event $H_2$ holds then the process $X$
	hits $\gamma_x$ before it exits  $D$. So we have
	$$ f(x''_1) \ge \E^{x_1''} ( f(X_{T_{\gamma_x} \wedge \tau_D}); H_2)  
	\ge p_1 \inf_{y \in \gamma_x} f(y) \ge p_1 C_1  f(x''_1) \ge 2 f(x''_1), $$
	a contradiction. (If $K \cap D_c = \emptyset$, we use the event $H_3$.)

	
	\medskip
	For the second part, let us first fix any point $x_0\in U_0$ and prove the argument for $f=G_D(x_0,\cdot)$. The argument for the case $f=G_D(x_0,\cdot)$ is similar; the main difference is in the
	definition of the path $\gamma$. In this case $f$ is harmonic in $D- \{x_0\}$ 
	(see \cite[Theorem 3.35]{Morter_Peres}), and so the path $\gamma_x$ on which we have
	$f(y) > C_1 f(x''_1)$ goes from $x$ to $x_0$. As $x_0 \in U_0$, there exists a path
	$\gam_0 \in D \cap B(0, \tfrac{1}{16})$ from $x_0$ to $(\tfrac{1}{32},0)$. 
	Let $S_3$ be the line segment between $(\tfrac{1}{32},0)$ and $(r_1,0)$.
	We then obtain a loop $\gamma$ in $D$ which contains the points $x, x_0, (\tfrac{1}{32},0)$,
	$(r_1,0)$, $z_1, x'_1,x$ by concatenating $\gamma_x$, $\gamma_0$, $S_3, S_1, L_1, \gamma'_x$.
	The construction of $\gamma$ gives that $\gamma \cap K= \emptyset$,
	and the remainder of the argument is the same as for the case $f=u$. 
\end{proof}

\begin{remark} \label{R:winding}
	{\rm 
		One can prove 
		formally  that the domains $D_a$ and $D_c$ are distinct using winding numbers. 
		We just give a sketch for the case of the function $f$. 
		Let $\gamma$ be the path between $x^*$ and $(1,0)$ constructed in the Lemma above,
		and $\gamma_1$ be a path in $\ol B(0,1)^c$ connecting $x^*$ and $(1,0)$. Let
		$\gamma_0$ be the closed path obtained by combining $\gamma$ and $\gamma_0$.
		
		Consider the contour integrals
		$ \oint_{\gamma_0} {dz}/{(z -z_b)} \hbox { for } b=a,c. $
		Let $\delta>0$  and assume that $|z_a- x_1''| = |z_c- x_1''| = \delta$. 
		As $z_a$ and $z_c$ are on opposite sides of $L_1$,
		if $\delta$ is small enough, then 
		$$  \Big | \int_{L_1}  \frac{dz}{(z -z_a)}  -   \int_{L_1}  \frac{dz}{(z -z_c)} \Big| \ge 1 . $$
		The path $\gamma'= \gamma_0  \setminus L_1$ is at least a distance $r_0/3$ from $x_1''$,
		and thus if $\delta$ is small enough the integrals
		$\int_{\gamma'} dz/(z-z_b)$ for $b=a,c$ will differ by less than $c \delta$.
		It follows that $z_a$ and $z_c$ are in different components of $\bR^2 \setminus \gamma_0$,
		and are therefore also  in different components of $B(0,1) \setminus \gamma$.
}\end{remark}

\medskip
The Lemma above controls the functions $u$ and $G_D(x, \cdot)$ (for $x\in U_0$) in the set $U_j$, $j=1,2$. 
We are ready to prove the main result of this section.

\medskip\noindent{\em{}Proof of Theorem 1.}
Let $u$, $v$ satisfy the hypotheses of the theorem.
Suppose first that $K \cap \pd B(0, \tfrac14)$ and $K \cap \pd B(0, \tfrac34)$
are non-empty.
Assume $x_i, z_i, \theta_i$ for $i=1,2$ and $w_0$ are as in Lemma \ref{rulce}. 

	{ \centering
\includegraphics[width=0.8\columnwidth]{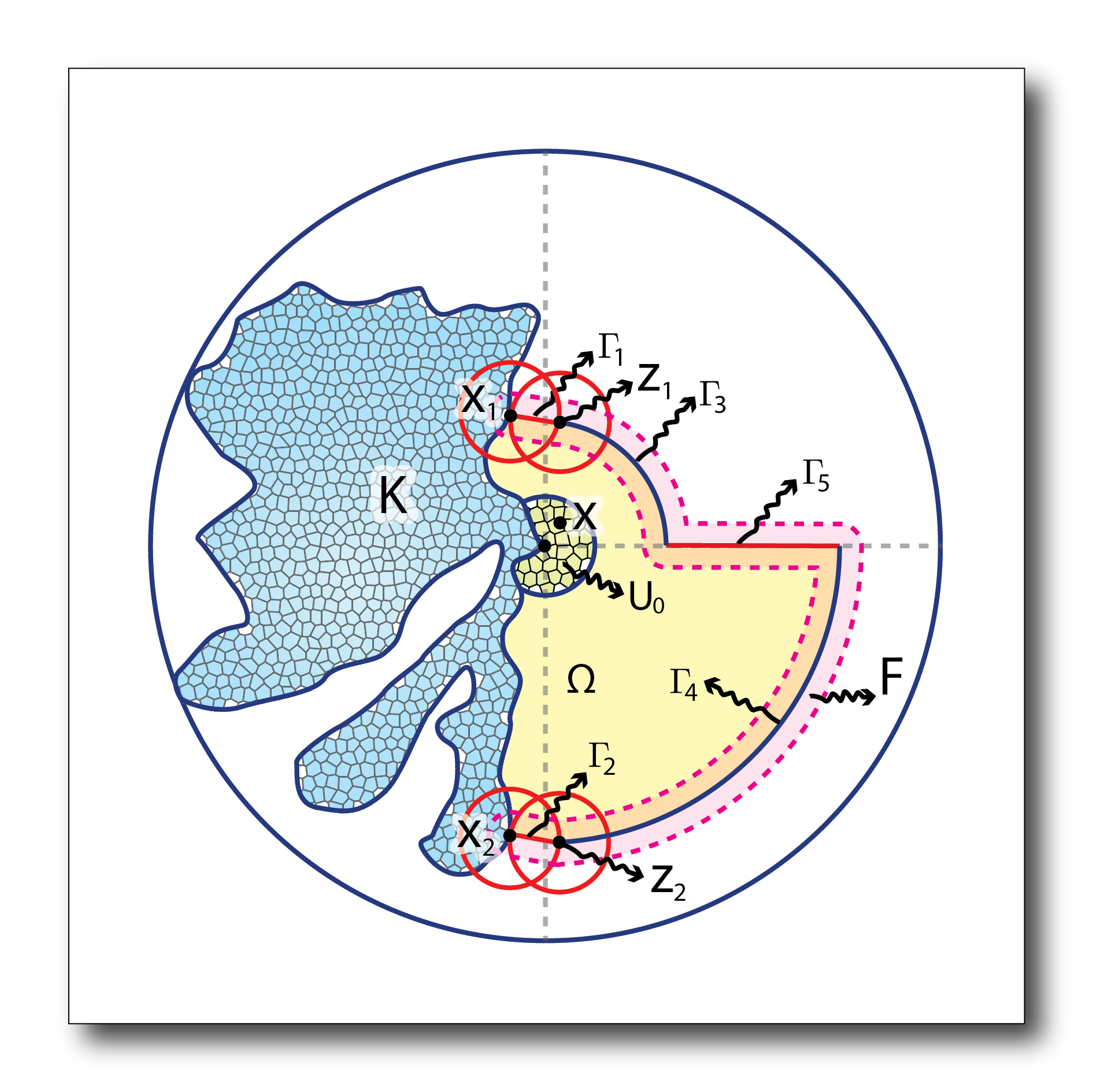}\\
\captionof{figure}{The sets $\Omega$ and $F$ }\label{fig4}
	}

\medskip
First we define the following 
paths which enclose the region we will work on:
\begin{align*}
	\Gamma_1&:=\{t\,x_1+(1-t) \,  z_1: t\in(0,1)\}, \\
	\Gamma_2&:=\{t\,x_2+(1-t) \,  z_2: t\in(0,1)\}, \\
	\Gamma_3&:=\{ r_1 e^{i\theta}: 0 \le \theta \le \theta_1\}, \\
	\Gamma_4&:=\{ r_2 e^{i\theta}: \theta_2 \le \theta \le  0 \}, \\
	\Gamma_5&=\{ (t,0), r_1 \le t \le r_2 \}.
\end{align*}
We write $\Gamma$ for the union of $\Gamma_1, \dots \Gamma_5$.			
Let $\Omega$ be the connected domain enclosed by the curves $\Gamma_j$, $j=1, \dots 5$
and the set $K$ and including the point $w_1=(\tfrac18,0)$. 

Let $x \in U_0\subset   \Omega$. 
Since $u$ is harmonic in $D$ and zero on $K$, 
\begin{align*}
	u(x)=\int_{\partial\Omega} u(y) \P^{x}({X_{\tau_{\Omega}}}\in dy)
	=\int_{\Gamma} u(y) \P^{x}({X_{\tau_{\Omega}}}\in dy).
\end{align*}
By the Carleson estimate Lemma \ref{rulce}, $u$ is bounded above by 
$c_1\, u(w_0)$ on $\Gamma_1\cup \Gamma_2$. By using a Harnack chain as in the 
proof of Lemma \ref{rulce} and the regular Harnack inequality, we also obtain that 
$u(y)\leq c_2 u(w_0)$ for $y\in \Gamma_3\cup \Gamma_4 \cup \Gamma_5$. 
Hence, by setting $c_3=\max\{c_1,c_2\}$,
\begin{align}\label{eqn_u}
	u(x) \le c_3 \,u(w_0)\, \P^{x}({X_{\tau_{\Omega}}}\in \Gamma).
\end{align}

Next, we define a tube $F$ around $\Gamma$ by
\bes  
F=\bigcup_{y\in \Gamma} B(y,\frac{r_0}{4}).
\ees 
By Urysohn's Lemma, there exists a smooth function $\psi$  with compact support in $F$ such that 
$\psi=1$ on $\Gamma$. 
We can choose $\psi$ so that $|\Delta \psi | \le C_1$.

For $x\in U_0$ we have

$$\int_{\partial \Omega} \psi(y)\,  \P^x({X_{\tau_{\Omega}}}\in dy)
	= \psi(x)+\int_{\Omega\cap \supp(\psi)} \Delta\psi(y) G_{\Omega}(x,y)\, dy.$$

(The functions on each side are harmonic in $\Omega$ and have the same 
boundary conditions -- see \cite[Lemma 1]{AikC} for details.)
This equation yields  
\begin{align}\label{eqn_prob_meas}
	\P^x({X_{\tau_{\Omega}}}\in \Gamma)
	&\leq \int_{\Omega\cap  F}  | \Delta\psi(y) | G_{D}(x,y)\, dy 
	\le C_2  \int_{\Omega\cap  F}  G_{D}(x,y)\, dy.
\end{align}

Now fix $x_0 \in U_0$. 	Thus
$G_{D}(x_0,\cdot)$ is a positive harmonic function in the domain 
$D - U_0$, which contains $\Omega\cap F$.  We have
\begin{align} \label{e:bdOmF}
	\partial (\Omega\cap F)
	&= [\partial (\Omega \cap F) \cap K] \cup \Gamma \cup [\Omega \cap \pd F ].
\end{align}
We now claim that 
\begin{align}\label{ge-1}
	G_{D}(x_0,y)\leq c_4 \,G_{D}(x_0,w_0) \hbox{ for } y \in \pd( \Omega \cap F ).
\end{align}
We consider in turn the three parts of the boundary given by \eqref{e:bdOmF}. 
If $y \in K$ then $G_{D}(x_0,y)=0$ so \eqref{ge-1} holds for any $x\in  \partial (\Omega \cap F) \cap K$. 
If $y \in \Gamma\cap B(x_j, r_0/2)$ then Lemma \ref{rulce} implies \eqref{ge-1}. If $y \in \Gamma - B(x_j, r_0/2)$ then y is at a positive distance from $K$ and hence \eqref{ge-1} holds by Harnack inequality applied on a chain of balls. So \eqref{ge-1} holds for any $y\in \Gamma$.

For the final part of the boundary, let $y \in \Omega \cap \pd F$. If $y \in B(x_1, r_0/2)$ then
$y \in U_1$ and we can we use Lemma \ref{rulce} again. A similar argument gives that \eqref{ge-1} holds if $y \in B(x_2, r_0/2)$. The remaining part of $\Omega \cap \pd F$
is a distance at least $c>0$ away from $K$, so 
using the Harnack inequality on a  Harnack chain we obtain \eqref{ge-1}, completing the
proof of the claim \eqref{ge-1}.

Once we have  \eqref{ge-1} the maximum principle gives that 
\be \label{e:ge-2} 
G_{D}(x_0,y)\leq c_4 \,G_{D}(x_0,w_0), \quad y \in F \cap \Omega.
\ee
Combining this with \eqref{eqn_prob_meas}, and using the fact that 
$|\Delta \psi| \le C_2$ on $F$ we obtain

\begin{align}\label{eqn_prob_meas_2}
	\P^{x_0}({X_{\tau_{\Omega}}}\in \Gamma )& 
	\leq c_6   \, G_{D}(x_0,w_0).
\end{align}
Combining \eqref{eqn_prob_meas_2} with \eqref{eqn_u} gives that
$$  u(x) \le c u(w_0) G_{D}(x,w_0) \hbox { for } x \in U_0. $$

For the final part of the proof, consider the circle $\partial B(w_0,r_0)$. By our assumption on the 
harmonic function $v$ and the Harnack inequality 
$$v(z)\geq c_7 \, v(w_0), \qquad z\in \partial B(w_0,r_0).$$
Moreover
$$G_{D}(z,w_0) \leq c_8, \qquad z\in \partial B(w_0,r_0)$$
and so
$$G_{D}(z,w_0) \leq (c_8/c_7) \, \frac{v(z)}{v(w_0)} , \qquad z\in \partial (D-B(w_0,r_0))$$
since $v$ is positive. By the maximum principle, the last inequality holds inside the domain 
$D-B(w_0,r_0)$ which includes $U_0$. Using this inequality, 
(\ref{eqn_u}) and (\ref{eqn_prob_meas_2}), we obtain
$$ u(x)\leq c_9 \, u(w_0)\, \frac{v(x)}{v(w_0)}\leq c_{10}\, u_n(w_0)\, \frac{v(x)}{v( w_0)} $$
where $c_9=c_3c_6c_8c'/c_7$ and we applied the Harnack inequality to have $v( w_0)\leq c'\,v(w_0)$.

Finally, if we switch the roles of $u$ and $v$ and the roles of $x$ and $y$ we also obtain
$$\frac{v(y)}{u(y)}\leq c_{10} \frac{v(w_0)}{u(w_0)}$$
which leads to the result
$$\frac{u(x)/v(x)}{u(y)/v(y)}\leq c_{10}^2= C_0.$$

Now suppose that $K \cap \pd B(0, \tfrac14)= K \cap \pd B(0, \tfrac34) = \emptyset$.
If $K \cap B(0,\tfrac14)= \emptyset$ then as $u,v$ are harmonic in $B(0,\tfrac14)$ the inequality
\eqref{e:bhp2b} follows from the Harnack inequality. So it remains to consider the
case when $K \subset B(0, \tfrac14)$. 

Let $\Gamma = \pd B(0, \tfrac12)$. Then for $x \in B(0, \tfrac{1}{16})-K$ we have 
$$ u(x) = \int_{\Gamma} u(y) \bP^x( X_{T_\Gamma \wedge \tau_D} \in dy ). $$
The Harnack inequality gives that
$$ C^{-1} u(w_0) \le u(y) \le C u(w_0) \hbox { for } y \in \Gamma, $$
and thus if $p(x) = \bP^x( T_\Gamma < \tau_D)$ we have
$$  C^{-1} u(w_0) p(x) \le u(x) \le  C u(w_0) p(x). $$
A similar inequality holds for $v$, and \eqref{e:bhp2b} follows immediately\qed

\begin{example} \label{E:dge3}
	{\rm
		The following example  shows that one cannot expect a similar 
		uniform BHP in higher dimensions.
		Let $d \ge 3$, $B=B(0,1)$,
		$K_0 = B \cap \bH_0$, where $\bH_0=\{x: \pi_1(x) =0\}$.
		Let $\delta$ be small and positive.
		Set 
		$$ K = K_0 - B(0, \delta), \quad D=B(0,1)-K. $$
		Thus $K$ is a $d-1$ dimensional plate with a small hole in the centre, and
		is connected. 
		Let $y$ be on the $x_1$ axis with $\pi_1(y)=1/4$.
		Let $u_-$ and $u_+$ be the harmonic functions in $D$ with boundary 
		condition 1 on $\partial B \cap \bH_-$ and  $\partial B \cap \bH_+$ 
		respectively, and zero boundary conditions elsewhere. Set
		$v= u_-+u_+$. So if  $\tau=\tau_D$ then we can write
		\begin{align*}
			u_-(x) = \bP^x( X_\tau \in \partial D \cap \bH_-, \tau< T_K),
			\quad v(x) = \bP^x( \tau < T_K). 
		\end{align*}
		By symmetry we have 
		$$ \frac{ u_-(0)}{v(0)} = 1/2. $$  
		
		On the other hand if $B'=B(0, \delta)$ then 
		$$ \bP^{y}( T_{B'} < \tau_D) \le  \bP^{y}( T_{B'} < \tau_B) \le c\delta^{d-2}. $$
		So we have
		$$ v(y) \asymp 1, \quad u_-(y) \le c \delta^{d-2}. $$
		Thus
		\be \label{e:bh-ce}
		\frac{ u_-(0)/v(0)}{u_-(y)/v(y)}  \ge c \delta^{2-d}. 
		\ee
		By continuity this inequality will also hold if $0$ is replaced by a point
		$x$ close to 0 with $\pi_1(x)>0$.
}\end{example}

\begin{example}  \label{E:Kconn}
	{\rm The same example, taking $K = \{ (y,0): \delta \le |y| <1 \}$, shows that
		one cannot drop the hypothesis that $K$ is connected from Theorem \ref{T:ubhp2}.
}\end{example}		

\begin{example}  \label{E:Ball}
	{\rm Now let $K = \{ (0,y) : |y| < 1- \varepsilon \}$ and $u^{(\varepsilon)}_\pm$ be as in Example \ref{E:dge3}.
		Let $r = \frac{1}{20}$ and $x_- =(-r,0)$, $x_+=(r,0)$. Then we have 
		$$ \lim_{\varepsilon \to 0} u^{(\varepsilon)}_-(x_-) =  \lim_{\varepsilon \to 0} u^{(\varepsilon)}_+(x_+)=p, 
		\quad
		\lim_{\varepsilon \to 0} u^{(\varepsilon)}_-(x_+) =  \lim_{\varepsilon \to 0} u^{(\varepsilon)}_+(x_-)=0, $$
		for some $p\in(0,1)$. Thus we cannot have a BHP which holds for all $x,y \in B(0, \tfrac{1}{16}) \cap D$.
}\end{example}					 
\newpage

\section{Uniform BHP for a harmonic function associated with cones in higher dimensions ($d\geq 2$)} 
\label{S:BHPuv}

In this section we will prove the Boundary Harnack Principle in $\R^d$ with $d \ge 2$ for two fundamental 
harmonic functions.

Let 
$$ \bH_- =\{ x=(x_1, \dots, x_d) \in \bR^d : x_1 < 0 \}, \quad  B^- = B(0,1) \cap \bH_-,$$
and define $\bH_+$, $B^+$ analogously. We write $B = B(0,1)$.
Let $K \subset \overline \bH_-$ be compact.
Set
$$ D = B(0,1) - K. $$
Let $\pi_1: \bR^d \to \bR^d$ be projection onto the $x_1$-axis, so
$\pi_1( (x_1, x_2, \dots , x_d) )=(x_1, 0, \dots, 0)$.
Let $W_\alpha$ be the cone 
$$ 	W_\alpha =\{z\in \Rd: |z-\pi_1(z) |<z_1 \tan(\alpha)\}. $$
Set
$$ W_\alpha(r) = B(0,r) \cap W_\al. $$

Write $\tau=\tau_D$ for the exit time of $X$ from $D$.

We define the functions
\begin{align} \label{e:uvdef}
	v(x) = \bP^x( X_\tau \in \pd D \cap K^c ), \quad
	u_\alpha(x) = \bP^x( X_\tau \in \pd D  \cap W_\alpha ). 
\end{align}
Thus $u_\al$ and $v$ are bounded, positive harmonic functions which vanish on $\pd K$, and have
boundary values 1 on $\pd D \cap K^c$ and $\pd D  \cap W_\alpha$ respectively. It is
clear that $u_\al \le v$ on $D$.
Both $u_\al$ and $v$ are bounded, positive and harmonic inside the domain $D$. 
Hence they satisfy the usual Harnack Inequality (see \cite[Theorem II.1.19]{bass})
in balls which are far enough from the boundary of $D$.
The main result of this section is  that these two functions satisfy 
a BHP with a constant which depends only on $d$ and $\alpha$.
(Note that since the geometry of the boundary of $K$ 
is not specified, classical results on the 
Boundary Harnack Principle such as \cite[Theorem III.1.2]{bass} do not apply.) \\

The main result of this section is  the following Theorem. 

\begin{theorem} \label{BHP}
	Let $\al \in (0,\pi/2)$.
	There is a constant $C=C(\alpha,d)>0$ depending only on $\alpha$ and $d$,
	and independent of $K$, such that 
	\be \label{e:bhp1}
	\frac{u_\al(x)/u_\al(y)}{v(x)/v(y)}\leq C
	\ee
	for any $x,y\in B(0,1/2) \cap \bH_+$.
\end{theorem}

\begin{remark} \label{R:alpha}
	{\rm 
		
		The usual Harnack inequality gives
		that $u_\al(\tfrac14,0,...,0) \ge C'(\al,d)$, and so, since $v \le 1$,  \eqref{e:bhp1} implies that
		\be  \label{e:bhp2}
		C''(\al,d) \le u_\al(x)/v(x) \le 1 \quad
		\ee
		for $x \in  B(0,1/2) \cap \bH_+$. On the other hand, given \eqref{e:bhp2} the inequality \eqref{e:bhp1} 
		follows immediately.
} \end{remark}

The proof of Theorem \ref{BHP} is based on  the Poison kernel of the unit ball which is harmonic inside the ball. First, we fix the point $y_0=(1,0,...,0)$ and define the Poisson kernel as
$$H(x):=\frac{1-|x|^2}{|x-y_0|^d}, \quad x\in B.$$
Before we prove the main theorem of this section, we state two short lemmas.

\begin{lemma}\label{comparison_H}
	For any $z\in  B^-$ and any point $x$ on the line segment from origin to the point $y_0$, the inequality
	$$H(z)\leq H(0) \leq H(x)$$
	holds. 
\end{lemma}

\begin{proof}
	Since $H(0)=1$ and $|z-y_0|\geq 1$ for any $z\in  B^-$, it is trivial to show $H(z)\leq H(0)$. Moreover, if $x$ is any point on the line segment from origin to the point $y_0$, then we have $1-|x|=|x-y_0|$. In this case,
	$$H(x)=\frac{1+|x|}{1-|x|} \, \frac{1}{|x-y_0|^{d-2}}$$
	where both terms are greater than or equal to 1 when $d\geq 2$. Hence the second inequality follows.
\end{proof}

For the next Lemma, we define the cross-section of the ball through the point $(\cos\alpha,0,...,0)$ as
$$S_\alpha=\{z\in B^+: |\pi_1(z)|=\cos\alpha\}.$$
\begin{lemma}\label{cross-section-comp}
	Assume $\alpha \in (0,\pi/2)$. For any $x\in S_\alpha$,	$$H(x)\leq\frac{1-\cos^2\alpha}{(1-\cos\alpha)^d}.$$
\end{lemma}
\begin{proof}
	By symmetry it is enough to prove the inequality for $x$ of the form $x=(\cos\alpha,t,0,...,0)$,
	where $t\in [0,\sin\alpha)$.
	
	For such points, the function $H$ becomes
	$$H(x)=\frac{1-\cos^2\alpha-t^2}{((1-\cos\alpha)^2+t^2)^{d/2}}.$$
	As a function of $t$, its derivative is negative for $t\in [0,\sin\alpha)$. Hence it attains its maximum value at $t=0$, which leads to the result. 
\end{proof}

Now, we are ready to prove the main result of this section.

\begin{proof}[Proof of Theorem \ref{BHP}]  
	By Remark 2 it is sufficient to show \eqref{e:bhp2} for $x \in B(0,1/2) \cap \bH_+$ . Since $u_\alpha$ is increasing with respect to $\alpha$, we may assume that $\alpha\in (0,\pi/3)$.
	
	We define the ball  $B^*$ obtained by shifting the center of the unit ball $B$ to  $(2\cos\alpha,0,...,0)$, and the domain $D^*_\alpha$ to be the intersection $B\cap B^*$. We note that the cross-section of the ball $B$ through the intersection of surfaces $\partial B^*$ and $\partial B$ is exactly $S_\alpha$	and $D^*_\alpha$ is symmetric about $S_\alpha$. 
	Moreover our choice of $\alpha$ ensures that $D^*_\alpha\subset B^+$.
	Let $L_\alpha$ be the line segment from the origin to the point $(\cos\alpha,0,...,0)$.
	
	By the symmetry of the domain $D^*_\alpha$ about $S_\alpha$, for any $x\in S_\alpha$, 
	\begin{align}\label{u_alpha_lowerbound}
		u_\alpha(x) \geq  \P^x( X_{\tau_{D^*_\alpha}}\in \pd D  \cap W_\alpha)=1/2.
	\end{align}
	
	Now we set 
	$$h_\alpha(x)=(1-v(x))+2\frac{1-\cos^2\alpha}{(1-\cos\alpha)^d}u_\alpha(x), \quad x\in \overline{D}.$$
	We claim that $h_\alpha(x)\geq H(x)$ on the set $(\partial D - W_\alpha)\cup S_\alpha$,
	which is the boundary of the domain  $\{x\in D: \pi_1(x)<\cos\alpha\}$.
	To see this, first we note $h_\alpha(x)=1\geq H(x)$ for $x\in \partial D \cap K$, by Lemma \ref{comparison_H}. Next, if $x \in \partial D-(K\cup W_\alpha)=\partial B-(K\cup W_\alpha)$ then $h_\alpha(x)=0= H(x)$. Finally, if $x\in S_\alpha$ then 
	$$h_\alpha(x)\geq 2\frac{1-\cos^2\alpha}{(1-\cos\alpha)^d}u_\alpha(x)\geq \frac{1-\cos^2\alpha}{(1-\cos\alpha)^d}\geq H(x)$$
	by Lemma \ref{cross-section-comp} and Inequality (\ref{u_alpha_lowerbound}). 
	
	Using the Maximum Principle for harmonic functions, we conclude that $h_\alpha(x)\geq H(x)$ also holds in 
	$\{x\in D: \pi_1(x)<\cos\alpha\}$, and so in particular for  $x\in L_\alpha$.
	
	Hence, by Lemma \ref{comparison_H}, for any $x\in L_\alpha$,
	\begin{align*}
		v(x)&\leq H(x)-1+v(x)\leq 2\frac{1-\cos^2\alpha}{(1-\cos\alpha)^d}u_\alpha(x),
	\end{align*}
	which leads to 
	\begin{align*}
		\frac{(1-\cos\alpha)^d}{2(1-\cos^2\alpha)} \leq 	\frac{u_\alpha(x)}{v(x)}\leq 1
	\end{align*}
	if we also recall that  $u_\alpha(x)\leq v(x)$. This proves \eqref{e:bhp2} for $x\in L_\alpha$ 
	for the case $\alpha\in(0,\pi/3)$. Since $u_\alpha$ is increasing with respect to $\alpha$, we have
	\begin{align}\label{e:esc1:3}
		c_{\alpha,d} \leq 	\frac{u_\alpha(x)}{v(x)}\leq 1  
	\end{align}
	for  $ x \in L_\alpha$ and  any $\alpha\in(0,\pi/2)$.

	\bigskip		
	Now let  $x \in B(0,1/2) \cap \bH_+$, and set $x'= x -\pi_1(x)$. 
	Let $W' = x' + W_\alpha$, 
	$A'=W'\cap\partial B(x',1/2)$, and write $\tau' = \tau_{B(x',1/2)}$.
	As the bound \eqref{e:esc1:3} does not depend on the set $K$, we can apply \eqref{e:esc1:3} to the ball $B(x’, 1/2)-K$ (with rescaling), to 
	deduce that
	\bes
	c'_{\alpha,d} \le  \P^x[ X_{\tau' } \in A' | \tau' < T_K] 
	= \frac{  \P^x[ X_{\tau' } \in A' ; \tau' < T_K] }{  \P^x[ \tau' < T_K] }.
	\ees
	Since $\tau' < \tau = \tau_{B}$ this implies
	$$   \P^x[ X_{\tau' } \in A' ; \tau' < T_K] \ge c'_{\alpha,d} v(x). $$
	Now by the standard Harnack inequality, 
	\bes
	\P^y [ X_\tau \in \pd D  \cap W_\alpha,  \tau < T_K] \ge c''_{\alpha,d}  \text { for } y \in A'. 
	\ees
	This can be seen by observing that $A'$ has a positive distance from the boundary of the region $B^+$ where this distance depends only on $\alpha$ and $d$ but not on the point $y$. Then
	\begin{align*}
		u_\alpha(x) &=  \P^x[ X_\tau \in \pd D  \cap W_\alpha , \tau < T_K]  \\
		&\ge   \P^x[ X_\tau \in \pd D  \cap W_\alpha , X_{\tau'} \in A', \tau' < T_K,  \tau < T_K]  \\
		&= \E^x \Big[ 1_{(  \tau' < T_K, X_{\tau'} \in A' )} \P^{X_{\tau'}}[ X_\tau \in \pd D  \cap W_\alpha,  \tau < T_K] \Big]\\
		&\ge c''_{\alpha,d} \P^x [ \tau' < T_K, X_{\tau'} \in A' ] \ge c''_{\alpha,d} c'_{\alpha,d} v(x). 
	\end{align*}
	Since we have $u_\alpha(x) \le v(x)$ everywhere, it follows that
	\bes
	C(\alpha,d)  \le \frac{u_\alpha(x)}{v(x)} \le 1 \text{ for } x \in \bH_+ \cap B(0,{1/2}).
	\ees
	This proves \eqref{e:bhp2}, and Theorem \ref{BHP} then
	follows from Remark \ref{R:alpha}.
	
\end{proof}

{\bf Acknowledgements}
	The first named author, Martin T. Barlow, was partially supported by NSERC (Canada). The second named author, Deniz Karl\i, was partially supported by NSERC (Canada) and partially by the BAP grant, numbered 20A101, at the I\c{s}\i k University, Istanbul, Turkey. We thank our referees for their comments, and in particular one referee
	for suggesting a considerable simplification of our proof of Theorem \ref{BHP}. We also
	thank P\i nar Karl\i\,Akg\"un for drawing the figures in this manuscript.

%
%



\end{document}